\documentclass[a4paper,11pt,reqno]{amsart} 
\usepackage{amsmath,amssymb,amsthm,latexsym}
\usepackage[numbers,sort&compress]{natbib}

\newcommand{\embed}{\hookrightarrow}  
\newcommand{\cembed}{\hookrightarrow^{\text{compact}}}
\newcommand{\charak}{\chi}  
  
\newcommand{\pfeil}{\mathop{\to}\limits}
\newcommand{\xpfeil}{\xrightarrow}  
\newcommand{\pfeillang}{\mathop{\longrightarrow}\limits}  
 
\newcommand{\where}{\,|\,}

\newcommand{\supm}{\mathop{\sup}\limits}

\newcommand{\schwach}{\rightharpoonup}  
\newcommand{\intg}{\mathop{\int}\limits}

\newcommand{\laplace}{\Delta}

\newcommand{\nz}{{\mathbb N}}

\newcommand{\rz}{{\mathbb R}}

\newcommand{\eps}{\varepsilon}  
\renewcommand{\phi}{\varphi}  

\theoremstyle{plain}
\newtheorem{theorem}{Theorem}[section]  
\newtheorem{corollary}[theorem]{Corollary}  
\newtheorem{lemma}[theorem]{Lemma}  

\newtheorem*{theorem*}{Theorem}

\theoremstyle{remark}  
\newtheorem{remark}[theorem]{Remark}

\newtheoremstyle{citing}
  {3pt}
  {3pt}
  {\itshape}
  {}
  {\bfseries}
  {.}
  {.5em}
  {\thmnote{#3}}

\theoremstyle{citing}
\newtheorem*{varthm}{}

\numberwithin{equation}{section}

\begin{document}
 
\title[Compact embeddings and indefinite semilinear elliptic problems]{Compact embeddings and\\
indefinite semilinear elliptic problems}

\author{Matthias Schneider}
\address{Fachbereich Mathematik\\  
        Johannes Gutenberg-Universit\"at\\  
        Staudinger Weg 9\\  
        55099 Mainz, Germany}  

\email{adi@mathematik.uni-mainz.de, schneid@sissa.it}

\date{March, 2001}  
\keywords{compact embeddings, concentration compactness, Emden Fowler, indefinite nonlinearities}
\subjclass{35J65, 35D05}
\begin{abstract}
Our purpose is to find positive solutions $u \in D^{1,2}(\rz^N)$ of the semilinear elliptic problem
$-\laplace u = h(x) u^{p-1}$ for $2<p$. The function $h$ may have an indefinite sign. Key ingredients are
a $h$-dependent concentration-compactness Lemma and a characterization of compact embeddings of 
$D^{1,2}(\rz^N)$ into weighted Lebesgue spaces.   
\end{abstract}

\maketitle

\section{Introduction}
We are interested in finding weak nonnegative solutions of Emden-Fowler type problems 
\begin{align}\label{eq:10}  
\begin{array}{l}  
-\laplace u - h(x) u^{p-1}=0 \qquad \text{ in } \rz^N,\\  
0 \lneqq u \in E:= D^{1,2}(\rz^N) \cap L^{p}(\rz^N,|h|).  
\end{array}  
\end{align}
From now on we make the assumption: 
\begin{align} \label{eq:13}
N \ge 3,\; p>2 \text{ and } h \in L^1_{loc}:\; h^+(x):= \max(0,h(x))\not \equiv 0.  
\end{align}
We denote by ${D^{1,2}(\rz^N)}$ the closure of ${C^{\infty}_c(\rz^N)}$ with respect to the norm 
${(\int |\nabla u|^{2})^{\frac{1}{2}}}$ in $L^{2^*}$. Moreover,
$L^p(\Omega,|h|)$ denotes the space of measurable functions $u$ satisfying
\[ \|u\|_{L^p(\Omega,|h|)}^p := \int_\Omega |h|\, |u|^p  = \|\, \charak_{\Omega}\,
|h|^{\frac{1}{p}}\, u\, \|_p^p< \infty.\]
$E$ is a Banach space equipped with the norm $\|u\|_E:= \|\nabla u\|_2 + \|u\|_{L^p(\rz^N,|h|)}$.\\
Furthermore, we assume $h$ to be symmetric with respect to some compact subgroup $G$ of $O(N)$, 
the group of orthogonal linear transformations in $\rz^{N}$, i.e.
\begin{align}
\label{eq:30}
(g_* h)(x):= h(g^{-1} x) = h(x) \quad \forall g \in G \quad \mbox{a.e. in }
\rz^{N}.  
\end{align}
We denote by $D^{1,2}_G(\rz^N)$ the subspace of $D^{1,2}(\rz^N)$ consisting of all $G$-symmetric 
functions and define $E_G:= D^{1,2}_G(\rz^N) \cap L^p(\rz^N,|h|)$.\\  
The basic requirements on the positive part of $h$, $h^+ := \max(0,h)$, are:
\begin{align}
\text{There is a } G-\text{symmetric } u \in C^\infty_c(\rz^N): \; \int h |u|^p > 0, \label{eq:11}\\
D^{1,2}_G(\rz^N) \text{ is continuously embedded in } L^p(\rz^N, h^+).\label{eq:12}  
\end{align}
There have been many studies of the equation in (\ref{eq:10}), mostly for radially symmetric nonnegative 
functions $h$. We shall mention among them the work of Ding and Ni \cite{DinNi85}, Gidas and Spruck
\cite{GidSpr81}, Kusano and Naito \cite{KusNai86}, Noussair and Swanson \cite{NouSwa89,NouSwa88,NouSwa90}.\\
Tshinanga \cite{Tshi96} (see also \cite{NouSwa88}) proved without any
symmetry assumptions for {\em nonnegative} functions $h$ the existence of a 
solution to (\ref{eq:10}) if
\begin{align}
\label{eq:14}
0 \not \equiv h(x) \le \frac{C}{(1+|x|^2)^a}, \; \frac{2N-2a}{N-2}<p<2^* := \frac{2N}{N-2} 
\text{ for some }C>0,\, 0<a<2.  
\end{align}
Noussair and Swanson \cite{NouSwa90} obtained a solution of (\ref{eq:10}) for {\em nonnegative} $h$ if
\begin{align}
\label{eq:15}
2<p<2^*,\; 0 \not \equiv h \in L^q \cap L^\infty: \; 1<q<\frac{2^*}{2^*-p}.  
\end{align}
Rother \cite{Rot90} solved (\ref{eq:10}) for sign changing, radially symmetric functions $h$ if
\begin{align}
\label{eq:16}
\begin{split}
&h \in L^1_{loc},\; 0 \neq h^+(|x|)= k_1(|x|) + k_2(|x|)\text{ for some } k_1, k_2 \in L^1_{loc},\\
&\exists f \in L^\infty:\; 0 \le k_1(|x|) \le f(x)
|x|^{\frac{(N-2)p-2N}{2}} \text{and } f(x) \xpfeil[x \to 0]{|x| \to \infty} 0,\\
&k_2(|x|) \text{ is nonnegative and } \int_0^\infty k_2(r) r^{N-1-p\frac{N-2}{2}}\, \text{dr} \; < \infty.  
\end{split}  
\end{align}
We generalize the above results to possibly sign-changing and non-radial functions $h$.
Weak solutions of (\ref{eq:10}) correspond to nonnegative critical points of the
associated energy functional $I\in C^1(E_G,\rz) \cap C^1(E,\rz)$ defined by 
\[I(u):= \frac{1}{2} \int |\nabla u|^2 - \frac{1}{p} \int h(x) |u|^p.\]
From (\ref{eq:11}) and (\ref{eq:12}) it may be concluded that
\[c := \inf_{\gamma \in \Gamma} \max_{t \in [0,1]} I(\gamma(t)) >0,
\text{ where } \Gamma :=\{\gamma \in C([0,1], E_G)\where   \gamma(0)=0,\,
I(\gamma(1))<0\}.\]
Thus the mountain pass Theorem provides a $(PS)_c$ sequence, i.e. a sequence $(u_n)_{n \in \nz}$ satisfying
\[I(u_n) \to c, \; I'(u_n) \to 0 \text{ as }n \to \infty.\]
We shall show that if $D^{1,2}_G(\rz^N) \embed L^p(\rz^N,h^+)$ is compact, then every $(PS)_c$
sequence contains a convergent subsequence. Consequently we have
\begin{varthm}[Theorem \ref{s4t1}]
Suppose (\ref{eq:13})-(\ref{eq:12}) hold. If $D^{1,2}_G(\rz^N) \embed L^p(\rz^N,h^+)$ is compact, then (\ref{eq:10}) has a 
nontrivial, nonnegative weak solution. 
\end{varthm}
Section \ref{sec:compactness} is devoted to the study of embeddings of $D^{1,2}(\rz^N)$ into
weighted Lebesgue spaces, e.g. we prove
\begin{varthm}[Theorem \ref{sec:compactness:t1}]
Suppose $k \in L^1_{loc}$ is a nonnegative function and $q>2$.
Then ${\displaystyle  D^{1,2}(\rz^{N})}$ is compactly embedded in ${\displaystyle  L^{q}(\rz^N,k)}$ if and only if 
the following three conditions are satisfied:
\begin{align*}  
{\supm_{\substack{x\in \rz^N\\0<\rho}} \rho^{(1-\frac{N}{2})q} \hspace{-0.5em}\intg_{B_\rho(x)} k \; <
\infty},\quad
{\supm_{\substack{x\in \rz^N\\0<\rho}} \rho^{(1-\frac{N}{2})q} \hspace{-1.5em}\intg_{B_\rho(x)\backslash B_R(0)} 
\hspace{-1.5em}k \quad \pfeillang^{R \to\infty} 0,}\quad 
{\supm_{\substack{x \in \rz^N\\0<\rho <\delta}} \rho^{(1-\frac{N}{2})q}\hspace{-0.5em} \intg_{B_\rho(x)} k
\quad \pfeillang^{\delta \to 0} 0}.  
\end{align*}  
\end{varthm}
Theorem \ref{sec:compactness:t1} and Theorem \ref{s4t1} generalize the above existence results for
(\ref{eq:10}) obtained in
\cite{NouSwa88,Rot90,Tshi96}, because (\ref{eq:14}), (\ref{eq:15}) and (\ref{eq:16}) are sufficient for 
the compactness of the inclusion of $D^{1,2}_G(\rz^N)$ in $L^p(\rz^N,h^+)$, as it is shown in 
Corollary \ref{sec:compactness:co1} and in \cite[Lem. 6]{Rot90}.\\
To deal with the non-compact case we follow the notation of Smets in \cite{Sme99}, 
where the linear case was considered, and define for nonnegative $k \in L^1_{loc}$, 
${\displaystyle x \in \rz^N}$ and $r>0$:  
\begin{align}
\label{eq:39}
\begin{split}
S_{r,k} &:= \inf \left\{ \|\nabla u\|_2^{2}\;\where  \; u \in D^{1,2}(\rz^N
\backslash B_r(0)), \, \int k |u|^{q} = 1 \right\}\\  
S^\infty_{k} &:= \sup_{r >0} S_{r,k} = \lim_{r \to \infty} S_{r,k} \\  
S^{x}_{r,k} &:= \inf \left\{ \|\nabla u\|_2^{2}\;\where \; u \in D^{1,2}(B_r(x)), 
\, \int k |u|^{q} = 1 \right\}\\  
S^{x}_{k} &:= \sup_{r >0} S^{x}_{r,k} = \lim_{r \to 0} S^{x}_{r,k}\\  
S^{*}_k &:= \inf_{x \in \rz^N} S_k^{x}.  
\end{split}
\end{align}  
If $D^{1,2}(\rz^N)$ is embedded in $L^q(\rz^N,k)$, all these quantities are bounded
away from zero, however, some may be infinite, e.g. we have
\begin{varthm}[Corollary \ref{sconcompl1}]
Suppose $k \in L^1_{loc}$ is nonnegative, $q>2$ and $D^{1,2}(\rz^N)$ is embedded in $L^q(\rz^N,k)$. Then 
\[S^{*}_k= S^{\infty}_{k} = \infty \text{  if and only if  }  D^{1,2}(\rz^N) \embed
L^{q}(\rz^N,k) \mbox{ is compact.}\]  
\end{varthm}
Let us introduce the compactness threshold $c_0$, defined by
\[ c_0:=\left(\frac{1}{2}-\frac{1}{p}\right) \inf_{x \in \rz^{N}\cup
\{\infty\}} \left\{|G_x| (S^{x}_{h^{+}})^{\frac{p}{p-2}} \right\}, \]  
where ${\displaystyle |G_x| = \#\{gx\where  g \in G\}}$ and ${\displaystyle
|G_\infty|:=1}$. With the help of a concentration compactness Lemma, given in Section \ref{sec:conc-comp-lemma},
it is possible to show that every $(PS)_c$-sequence contains a convergent subsequence if $c<c_0$. This
is done in Section \ref{palais} and leads to
\begin{varthm}[Theorem \ref{s4t2}]
Suppose $D^{1,2}(\rz^N)$ is continuously embedded in $L^p(\rz^N,h^+)$ 
and there is an ${\displaystyle u \in E_G}$ such that  
\[\int h |u|^{p}>0 \mbox{ and } \max_{0 \le t< \infty} I(tu) \le c_0.\]  
Then (\ref{eq:10}) is solvable. 
\end{varthm}
Our approach is related to the work of Bianchi, Chabrowski and Szulkin \cite{BiaChaSzu95}, where the
case $p=2^*$ and $h \in L^\infty$ was considered. Our results for $p=2^*$ are slight improvements of \cite{BiaChaSzu95}, 
because in our setting $h^-$ does not need to be bounded.\\
Section \ref{exa} presents some examples illustrating our results, e.g.
\begin{varthm}[Corollary \ref{sec:examples:c4}]
Consider the equation
\begin{align}
\label{eq:40}
-\laplace u = (1+|x|)^{-\delta} |u|^{p-2}u, \quad 0\neq u \in D^{1,2}(\rz^N) \cap
 L^p(\rz^N,(1+|x|)^{-\delta}).
\end{align}
\begin{itemize}
\item[(i)] (\ref{eq:40}) has no solution $u \in C^2(\rz^N)$ if $2<p<2^*$ and $\delta \le N
  -\frac{p}{2}(N-2)$.
\item[(ii)] (\ref{eq:40}) has infinitely many $C^2(\rz^N)$-solutions if $2<p<2^*$ and $\delta > N
  -\frac{p}{2}(N-2)$. At least one solution is strictly positive in $\rz^N$.  
\end{itemize}
\end{varthm}

\section{Compactness}
\label{sec:compactness}  
${\displaystyle D^{1,2}(\rz^{N})}$ is embedded in ${\displaystyle L^{2^{*}}(\rz^{N})}$ but
not in ${\displaystyle L^{q}(\rz^{N})}$ for any other $q$. However we have an
embedding if we replace ${\displaystyle L^{q}(\rz^{N})}$ by a suitable
weighted Lebesgue space. Results concerning existence or compactness 
of such embeddings are obtained by Mazja \cite{Maz85}, Adams \cite{Ada73}, Berger and
Schechter \cite{BerSche72}.
For $q>2$ it is known (\cite[1.4.1]{Maz85} or \cite{Ada73}), that
there are positive constants $c_1(N,q)$, $c_2(N,q)$ such that  
\begin{equation}  
c_1 \supm_{u \in D^{1,2}\backslash\{0\}} \frac{\|k^{\frac{1}{q}} u\|_q}{\|\nabla
u\|_2}\le \supm_{\substack{x \in \rz^N\\0<\rho}} \rho^{(1-\frac{N}{2})} \left( \intg_{B_\rho(x)} k
\right)^{\frac{1}{q}} \le c_2 \supm_{u \in D^{1,2}\backslash\{0\}}
\frac{\|k^{\frac{1}{q}} u\|_q}{\|\nabla u\|_2}.  
\label{mja1}  
\end{equation}  
It is also shown in \cite[1.4.6]{Maz85} for $q>2$: 
\begin{equation}  
W^{1,2}(B_R(0)) \to L^{q}(B_R(0),k) \quad \mbox{is compact if} \quad  
\supm_{\substack{x \in \rz^N\\0<\rho <\delta}} \rho^{(1-\frac{N}{2})q}
\hspace{-1em}\intg_{B_\rho(x)\cap B_R(0)}\hspace{-1em} k \quad
\pfeillang^{\delta \to 0} 0.  
\label{mja2}  
\end{equation}
We give a version for all of $\rz^N$.   
\begin{theorem}
\label{sec:compactness:t1}
Let $q>2$ and $k$ a nonnegative, measurable function. Then   
${\displaystyle  D^{1,2}(\rz^{N})}$ is compactly embedded in ${\displaystyle  L^{q}(\rz^N,k)}$ if and only if 
the following three conditions are satisfied:
\begin{align}
\label{eq:17}  
\supm_{\substack{x\in \rz^N\\0<\rho}} \rho^{(1-\frac{N}{2})q} \intg_{B_\rho(x)} k \; <
\infty,\\  
\label{eq:18}
\supm_{\substack{x\in \rz^N\\0<\rho}} \rho^{(1-\frac{N}{2})q} \hspace{-1em}
\intg_{B_\rho(x)\backslash B_R(0)} \hspace{-1em} k \quad \pfeillang^{R \to
\infty} 0,\\ 
\label{eq:19} 
\supm_{\substack{x \in \rz^N\\0<\rho <\delta}} \rho^{(1-\frac{N}{2})q} \intg_{B_\rho(x)} k
\quad \pfeillang^{\delta \to 0} 0.  
\end{align}  
\end{theorem}  
\begin{proof}[Proof (sufficiency)]  
Using (\ref{mja1}) and (\ref{eq:17}), we see that $D^{1,2}(\rz^{N})$ 
is continuously embedded in $L^{q}(\rz^{N},k)$. (\ref{eq:19})
together with (\ref{mja2}) show  
\[ D^{1,2}(B_R(0)) \embed W^{1,2}(B_R(0)) \cembed L^{q}(B_R(0),k)
\pfeil^{ext} L^{q}(\rz^{N},k)\]  
with compactness in the middle and where $ext$ is understood by extending the
function by zero outside ${\displaystyle B_R(0)}$.\\  
Consider the operator $I_R$ defined by 
${\displaystyle I_R(u) := u \cdot \eta(\frac{x}{R})}$, where $R>0$,
$\eta \in C^{1}(\rz^{N},[0,1])$ with ${\displaystyle \eta|_{B_1(0)} \equiv
1}$ and compact support in ${\displaystyle B_2(0)}$. Then  
\[ D^{1,2}(\rz^{N}) \pfeillang^{I_R} D^{1,2}(B_{2R}(0)) \cembed
L^{q}(\rz^{N},k)\]  
is compact. We use again (\ref{mja1}) to get  
\begin{align*}  
\|k^{\frac{1}{q}}(u-I_R(u))\|_q &\le  \|k^{\frac{1}{q}} \charak_{\rz^{N}\backslash B_R(0)} u\|_q \\  
&\le c_1^{-1} \left(\supm_{\substack{x \in \rz^N\\0<\rho}}  \rho^{(1-\frac{N}{2})q} \hspace{-1em}
\intg_{B_\rho(x)\backslash B_R(0)} k \right)^{\frac{1}{q}} \|\nabla u\|_2 =
o(1) \|\nabla u\|_2  
\end{align*}  
as $R \to \infty$. Thus the inclusion of $D^{1,2}(\rz^N)$ in $L^q(\rz^N,k)$ is compact 
as a limit of compact operators.
\end{proof} 
\begin{proof}[Proof (necessity)]
Because the inclusion of $D^{1,2}(\rz^N)$ in $L^q(\rz^N,k)$ is bounded (\ref{eq:17}) holds.
Let ${\displaystyle D:= \{u \in D^{1,2}(\rz^{N}) \where \|\nabla u\|_2\le 1\} }$. The set
$D$ is relatively compact in ${\displaystyle
L^{q}(\rz^{N},k)}$. Let $\eps>0$. The relative compactness gives rise to  
\begin{eqnarray}  
\exists R>0 : && \intg_{\rz^{N}\backslash B_{R}(0)} k \cdot |u|^{q} \le \eps
\quad \forall u \in D \label{iR}\\  
\exists \delta>0 : && \intg_{B_{2\delta}(x)} k \cdot |u|^{q} \le \eps \quad
\forall u \in D, \quad \forall x \in \rz^{N}. \label{iiR}  
\end{eqnarray}  
We show only (\ref{iiR}): We suppose the contrary and get a sequence
$(x_n,u_n)_{n \in \nz}$ such that for all $n \in \nz$ there holds  
${\int_{B_{\frac{1}{n}}(x_n)} k |u_n|^{q} > \eps}$. (\ref{iR})
implies $(x_n) \subset B_{R}(0)$ and we may assume $x_n \to x_0$ and
$u_n \to u_0$ in $L^{q}(\rz^N,k)$. Consequently we have  for all sufficiently large $n$ 
${\int_{B_{\frac{1}{n}}(x_0)} k |u_0|^{q} > \frac{\eps}{2}}$; a contradiction.\\  
Let $\Omega \subset \rz^{N}$ and $\rho >0$. Consider ${\displaystyle
u_{\rho,x}(y):= \eta(\frac{y-x}{\rho})}$. We obtain  
\begin{eqnarray*}  
\intg_{B_\rho(x)\backslash \Omega}k &\le& \intg_{B_\rho(x)\backslash \Omega}k
|u_{\rho,x}|^{q}\le \intg_{B_\rho(x)\backslash \Omega} k
\frac{|u_{\rho,x}|^{q}}{\|\nabla u_{\rho,x}\|_2^{q}} \, \|\nabla
u_{\rho,x}\|_2^{q} \\  
&=& \intg_{B_\rho(x)\backslash \Omega} k \frac{|u_{\rho,x}|^{q}}{\|\nabla
u_{\rho,x}\|_2^{q}} \, \rho^{-q} \rho^{\frac{N}{2}q} \|\nabla \eta\|_2^{q}  
\end{eqnarray*}  
and observe that ${\displaystyle \frac{u_{\rho,x}}{\|\nabla u_{\rho,x}\|_2}
\in D}$. If we take $B_R(0) \subset \Omega$, we get (\ref{eq:18}) and if we take
$\Omega= \emptyset$ and $0<\rho<\delta$, we get (\ref{eq:19}).  
\end{proof}
As an easy consequence of Theorem \ref{sec:compactness:t1} (see \cite{SchnDis01}) we obtain
\begin{corollary} \label{sec:compactness:co1}  
Suppose $2<q<2^*$ and $k \in L^1_{loc}$ is a nonnegative function.\\ 
Then $D^{1,2}(\rz^N) \embed L^q(\rz^N,k)$ is compact if one of the following conditions is satisfied
\begin{align}
&\int k^{\frac{2^{*}}{2^{*}-q}}< \infty, \label{eq:20}\\
&\exists f \in L^\infty: \; f(x) \xpfeil[x\to 0]{|x| \to \infty} 0 \text{ and }
k(x) \le f(x) |x|^{\frac{(N-2)q-2N}{2}}. \label{eq:21}
\end{align}
\end{corollary}
The following Theorem gives some conditions ensuring that $D^{1,2}(\rz^N)\cap L^p(\rz^N,h)$ is
compactly embedded in $L^q(\rz^N,k)$. Furthermore, it leads to sufficient conditions for $D^{1,2}(\rz^N)$ to 
be embedded in $L^q(\rz^N,k)$ for $1\le q\le 2$ (see Corollary \ref{sec:compactness:co2} below).
\begin{theorem}
\label{sec:Compactness:t2}  
Suppose $\min(p,2^{*})>q\ge 1$ and $h,k$ are nonnegative measurable functions which satisfy   
\begin{align*}  
&\exists R >0: \quad h>0 \; \mbox{almost everywhere in } \; \Omega_k \backslash B_R(0),\\  
&k\in L_{loc}^{\frac{2^*}{2^*-q}}(\rz^N) \text{ and }\intg_{\Omega_k \backslash B_R(0)} \left[k
\left(\frac{k}{h}\right)^{\frac{q-z}{p-q}} \right]^{\frac{2^{*}}{2^{*}-z}} \;
< \infty  
\end{align*}
for some $0\le z \le q$, where ${\displaystyle \Omega_k := \{x \in \rz^N \where k(x) \neq 0\}}$.\\
Then $D^{1,2}(\rz^{N})\cap L^{p}(\rz^N,h)$  is compactly embedded in $L^{q}(\rz^N,k)$.  
\end{theorem}  
\begin{proof}  
We follow \cite[Lem 2.3]{AlaTar96} and use the elementary upper bound for
$r>s\ge0$, ${\displaystyle u \in \rz}$, $k\ge 0$ and $h>0$:  
\begin{equation}  
k |u|^{s} - h|u|^{r} \le C(r,s) k \left(\frac{k}{h}\right)^{\frac{s}{r-s}}.  
\label{s1p2f1}  
\end{equation}
We use Theorem \ref{sec:compactness:t1} to see that $D^{1,2}(B_R(0))$ is compactly embedded in $L^q(\rz^N,k)$. 
If $q>2$, we may factorize the inclusion as follows 
\[D^{1,2}(B_R(0)) \embed D^{1,2}(\rz^N) \cembed L^{q}(\rz^{N},k \charak_{B_R(0)}) 
\xpfeil{\cdot \charak_{B_R(0)}} L^{q}(\rz^{N},k).\]  
If $q\le 2$, then we fix $q_1$ between $2$ and $2^*$ and notice that 
\[{k_1}:= k^{\frac{2^*-q_1}{2^*-q}}\cdot \charak_{B_R(0)} \in L^{\frac{2^*}{2^*-q_1}}.\] 
Hence we may write the inclusion of $D^{1,2}(B_R(0))$ in $L^q(\rz^N,k)$ as  
\[D^{1,2}(B_R(0)) \embed D^{1,2}(\rz^N) \cembed  L^{q_1}(\rz^N, k_1) 
\xpfeil{\cdot \charak_{B_R(0)}} L^q(\rz^N,k).\] 
The last multiplication operator is bounded due to H\"older's inequality.\\
Using again ${\displaystyle I_R}$ as in the proof of Theorem \ref{sec:compactness:t1}, we see  
\[D^{1,2}(\rz^{N}) \cap L^{p}(\rz^{N},h) \pfeillang^{I_R} D^{1,2}(B_{2R}(0))\cembed L^{q}(\rz^{N},k)\]  
is also compact.\\  
The bound in (\ref{s1p2f1}) allows to calculate for all $\epsilon>0$ and ${\displaystyle u \in
D^{1,2}(\rz^{N}) \cap L^{p}(\rz^{N},h)}$ with 
${\displaystyle (\|\nabla u\|_2^{2}+ \|h^{\frac{1}{p}}u\|_p^{2})^{\frac{1}{2}}} =1$:  
\begin{align*}  
\int k |u -I_R(u)|^q &\le \intg_{\rz^{N}\backslash
B_R(0)}\hspace{-1ex} k |u|^{q} = \intg_{\rz^{N}\backslash B_R(0)}\hspace{-1ex} 
(k |u|^{q} - \epsilon h |u|^{p}) + \epsilon \intg_{\rz^{N}\backslash B_R(0)}\hspace{-1ex} h |u|^{p} \\  
&\le \intg_{\rz^{N}\backslash B_R(0)}\hspace{-1ex} (k |u|^{q-z} - \epsilon h |u|^{p-z})
|u|^{z} + \epsilon \intg_{\rz^{N}\backslash B_R(0)}\hspace{-1ex} h |u|^{p} \\  
&\le \epsilon^{-\frac{q-z}{p-q}}C(p,q,z) \intg_{\Omega_k\backslash B_R(0)} \left[k
\left(\frac{k}{h}\right)^{\frac{q-z}{p-q}} \right] |u|^{z} + \epsilon
\intg_{\rz^{N}\backslash B_R(0)} h |u|^{p}\\  
&\le \epsilon^{-\frac{q-z}{p-q}} C(p,q,z) \left(\intg_{\Omega_k\backslash
B_R(0)} \left[k \left(\frac{k}{h}\right)^{\frac{q-z}{p-q}}
\right]^{\frac{2^{*}}{2^{*}-z}}\right)^{\frac{2^{*}-z}{2^{*}}}\hspace{-3ex}
\|u\|_{2^{*}}^{z} + \epsilon \hspace{-3ex}\intg_{\rz^{N}\backslash B_R(0)}\hspace{-2ex} h |u|^{p}\\  
&\le \epsilon^{-\frac{q-z}{p-q}} C(p,q,z) \left(\int_{\Omega_k\backslash
B_R(0)} \left[k \left(\frac{k}{h}\right)^{\frac{q-z}{p-q}}
\right]^{\frac{2^{*}}{2^{*}-z}}\right)^{\frac{2^{*}-z}{2^{*}}} + \epsilon.  
\end{align*}  
The integral term tends to zero for $R \to \infty$.\\
Hence the inclusion of $D^{1,2}(\rz^{N})\cap L^{p}(\rz^{N},h)$
in $L^{q}(\rz^{N})$ is compact as a limit of compact operators.
\end{proof}  

\begin{corollary}
\label{sec:compactness:co2}
Suppose $1\le q< 2^*$ and $k \in L^{\frac{2^*}{2^*-q}}_{loc}$ is a nonnegative function. Then
$D^{1,2}(\rz^N)$ is compactly embedded in $L^q(\rz^N,k)$ under the following
condition
\begin{align}
1 \le q < p:=\frac{2(N-\delta)}{N-2} \text{ and } \int k^{\frac{p}{p-q}} |x|^{\frac{\delta q}{p-q}} < \infty
\text{ for some }\delta: \;0 \le \delta \le 2. \label{eq:24}     
\end{align}
\end{corollary}
\begin{proof}
$D^{1,2}(\rz^N)$ is continuously embedded in $L^p(\rz^N,|x|^{-\delta})$ if (see (\ref{eq:4}) below) 
\[0 \le \delta \le 2 \text{ and } p= \frac{2(N-\delta)}{N-2}.\] 
Consequently (\ref{eq:24}) and Theorem \ref{sec:Compactness:t2} with $h(x):= |x|^{-\delta}$ and $z=0$ implies
\[D^{1,2}(\rz^N)= D^{1,2}(\rz^N) \cap L^p(\rz^N,|x|^{-\delta}) \cembed  L^q(\rz^N,k).\]  
\end{proof}
\section{A concentration compactness Lemma}
\label{sec:conc-comp-lemma}  
The following Lemma, which is closely related to 
\cite[Lem. 2.1]{Sme99} and \cite[Lem. I.1]{Lions85}, analyses
the possible non-compactness of an embedding of $D^{1,2}(\rz^N)$ in $L^q(\rz^N,k)$ in terms of the
quantities $S^{x}_k$ and $S^{\infty}_{k}$ defined in (\ref{eq:39}).     
\begin{lemma} \label{sec:conc-comp-lemma-l2} 
Suppose $q>2$ and $D^{1,2}(\rz^N)$ is continuously embedded in $L^q(\rz^N,k)$ for
some nonnegative $k \in L^1_{loc}$. Furthermore, let $(u_n)_{n \in \nz}$ be bounded in 
${\displaystyle D^{1,2}(\rz^N)}$. Up
to a subsequence we may assume: ${\displaystyle u_n \schwach u}$ weakly in
${\displaystyle D^{1,2}(\rz^N)}$ and additionally ${\displaystyle |\nabla u_n- \nabla u|^{2}
\schwach \tilde\mu}$, ${\displaystyle |\nabla u_n|^{2} \schwach \mu}$,
${\displaystyle k |u_n|^{q} \schwach \nu}$ and ${\displaystyle k |u_n- u|^{q}
\schwach \tilde\nu}$ weakly in the sense of measures, where $\mu, \tilde\mu$ and
$\nu$ are bounded nonnegative measures. Define  
\begin{eqnarray*}  
\mu_\infty &:=& \lim_{R \to \infty} \limsup_{n \to \infty} \int_{|x|>R} |\nabla u_n|^{2} \\  
\nu_\infty &:=& \lim_{R \to \infty} \limsup_{n \to \infty} \int_{|x|>R} k |u_n|^{q}.  
\end{eqnarray*}  
Then  
\begin{itemize}  
\item[(1)] ${\displaystyle \mu_{\infty} \ge S^{\infty}_{k}\, \nu_\infty^{2/q}}$,  
\item[(2)] There exists an at most countable set $J$, a family
${\displaystyle \{x_j\where  j \in J\}}$ of distinct points in ${\displaystyle
\rz^{N}}$ and a family ${\displaystyle \{\nu_j\where  j \in J\}}$ of positive
numbers such that  
\[ \nu = k |u|^{q} \mbox{ dx} + \sum_{j \in J} \nu_j \delta_{x_j} \]  
where $\delta_{x}$ is the Dirac measure of mass $1$ concentrated at $x \in \rz^{N}$,  
\item[(3)] There holds  
\[ \mu \ge |\nabla u|^{2} \mbox{ dx} + \sum_{j \in J} \mu_j \delta_{x_j},\]  
where ${\displaystyle \mu_j \ge S_k^{x_j}\, \nu_j^{2/q}}$ for all $j \in J$,  
\item[(4)] ${\displaystyle \limsup_{n \to \infty} \|k^{1/q} u_n\|_q^{q} = \|k^{1/q}
u\|_q^{q} + \sum_{j \in J} \nu_j + \nu_\infty }$.  
\end{itemize}  
\end{lemma}  
\begin{proof}  
Let $\{x_j\where  j \in J\}$ be the atoms of $\tilde \nu$ and decompose
${\displaystyle \tilde\nu = \nu_0 + \sum_{j \in J} \nu_j \delta_{x_j}}$, where
$\nu_0$ is nonnegative and free of atoms. 
Because ${\displaystyle \int d\tilde\nu < \infty}$, $J$ is
at most countable. For each ${\displaystyle x \in \{x_j\where j \in J\}}$ there is
a sequence $(r_l)_{l \in \nz}$ of positive numbers converging to zero such
that  
\[S_{r_l,k}^{x} \ge  
\left\{  
\begin{array}{lr}  
S^{x}_{k} - \frac{1}{l} & S^{x}_{k}<\infty\\  
l & S^{x}_{k} = \infty  
\end{array}  
\right.. 
\]  
Let ${\displaystyle 0\le \psi_l \in C^{\infty}_c(B_{r_l}(x))}$ with
${\displaystyle \|\psi_l\|_\infty = 1 = \psi_l(x)}$, then  
\begin{align*}  
\tilde\mu(\{x\}) &= \lim_{l \to \infty} \tilde\mu(\psi_l^{2}) = \lim_{l \to
\infty} \lim_{n \to \infty} \int |\nabla(u_n-u)|^{2} \psi_l^{2}\\  
&= \lim_{l \to \infty} \lim_{n \to \infty} \int |\nabla((u_n-u)\psi_l)|^{2}
 \tag{$\text{because }u_n \to u \text{ in } L^{2}_{loc}$}\\   
&\ge \lim_{l \to \infty} \left\{ S^{x}_{r_l,k} \limsup_{n \to \infty} \left(
\int k |u_n -u|^{q} \psi_l^{q} \right)^{2/q} \right\}\\  
&= \lim_{l \to \infty} S^{x}_{r_l,k} \tilde\nu(\psi_l^{q})^{2/q} = S^{x}_k
\tilde\nu(\{x\})^{2/q}.  
\end{align*}  
The above calculation also shows that  
\begin{equation}  
\left(\int |\psi|^{q} d\tilde\nu\right)^{2/q} \le C \int |\psi|^{2}
d\tilde\mu \qquad \forall \psi \in C^{\infty}_c(\rz^N).  
\label{4.11}  
\end{equation}  
This implies that ${\displaystyle \nu_0}$ is absolutely continuous with
respect to $\tilde\mu$. By the Radon-Nikodym Theorem there is a nonnegative $f \in
L^{1}(\rz^N , d\tilde\mu)$ such that ${\displaystyle d\nu_0 = f d\tilde\mu}$
and for $\tilde\mu$-almost every $x \in \rz^N$  
\[ f(x) = \lim_{r \to 0} \left(\frac{\nu_0(B_r(x))}{\tilde\mu(B_r(x))}\right).\]  
If $x$ is not an atom of $\tilde\mu$, we use (\ref{4.11}) to get  
\[f(x)^{2/q} = \lim_{r \to 0}
\left(\frac{\nu_0(B_r(x))^{2/q}}{\tilde\mu(B_r(x))^{2/q}}\right)  
\le C \lim_{r \to 0} \tilde\mu(B_r(x))^{\frac{q-2}{q}} =0.\]  
Because the atoms of $\tilde\mu$ are at most countable and $\nu_0$ has no
atoms, we see $\nu_0 =0$.\\  
We use the inequality ${\displaystyle |(a-b)^{2}-a^{2}| \le \epsilon
(a-b)^{2}+ c(\epsilon) b^{2}}$ to derive  
\begin{align*}  
\left|\int |\nabla(u_n-u)|^{2} \psi_l^{2} \right.&- \left.\int |\nabla
u_n|^{2}\psi_l^{2}\, \right|\\  
&\le \epsilon \int |\nabla(u_n-u)|^{2} \psi_l^{2} + c(\epsilon) \int |\nabla
u|^{2}\psi_l^{2}\\  
&\le \epsilon C + c(\epsilon) o(1)_{l \to \infty}.
\end{align*}  
Letting $l \to \infty$ we get that 
${\displaystyle \tilde\mu(\{x\}) = \mu(\{x\})}$. Because of the
weak lower semi-continuity we have $\mu \ge |\nabla u|^{2} dx$. Finally the
Brezis-Lieb Lemma \cite{BreLie83} implies  
\[k (|u_n|^{q}-|u|^{q}) dx =  k (|u_n - u|^{q}) dx + o(1)_{n \to \infty}.\]  
Thus claims (2) and (3) are proved.\\
Let $R>0$ and $\psi_R \in C^\infty(\rz^N)$ such that $\psi_R \equiv 0$ in $B_R(0)$, 
$\psi_R \equiv 1$ in $\rz^N\backslash B_{R+1}(0)$ and $0 \le \psi_R\le 1$ everywhere.\\
Because $u_n \to u$ in $L^2_{loc}$ we have
\begin{align*}
\mu_\infty &= \lim_{R \to \infty} \lim_{n \to \infty} \int |\nabla u_n|^2 \psi_R^2
= \lim_{R \to \infty} \lim_{n \to \infty} \int |\nabla(u_n \psi_R)|^2\\
&\ge  \lim_{R \to \infty} S_{R,k} \lim_{n \to \infty} \left(\int k |u_n \psi_R|^q\right)^{2/q}
= S_{\infty, k} \, \nu_\infty^{2/q}.
\end{align*}
To show (4) we use again the Brezis-Lieb Lemma and get
\begin{align*}
\lim_{R \to \infty} \limsup_{n \to \infty} \int k |u_n|^q (1\!\!-\!\!\psi_R) &=
\lim_{R \to \infty} \limsup_{n \to \infty}\left( \int k |u_n\!\!-\!\!u|^q (1\!\!-\!\!\psi_R) + 
\int k |u|^q (1\!\!-\!\!\psi_R)\right)\\
&= \sum_{j \in J} \nu_j + \int k |u|^q.     
\end{align*}
Finally we deduce 
\begin{align*}
\limsup_{n \to \infty} \int k |u_n|^q &= \lim_{R \to \infty}\limsup_{n \to \infty} 
\left(\int k |u_n|^q (1-\psi_R) + \int k |u_n|^q \psi_R \right) \\
&= \sum_{j \in J} \nu_j + \int k |u|^q + \nu_\infty.
\end{align*}
\end{proof}
\begin{corollary}
\label{sconcompl1}  
Suppose $q>2$ and  $k\in L^1_{loc}$ nonnegative such that ${\displaystyle D^{1,2}(\rz^N)}$ 
is continuously embedded in ${\displaystyle L^{q}(\rz^N,k)}$. Then  
\[S^{*}_k= S^{\infty}_{k} = \infty \text{  if and only if  }  D^{1,2}(\rz^N) \embed
L^{q}(\rz^N,k) \mbox{ is compact.}\]  
\end{corollary}
\begin{proof}[Proof (sufficiency)]
Suppose $(u_n)_{n \in \nz}$ is bounded in $D^{1,2}(\rz^N)$ such that $u_n \schwach 0$. Since $S^{*}_k=
S^{\infty}_{k} = \infty$, Lemma \ref{sconcompl1} shows that $\int k |u_n|^q \to 0$ as $n \to \infty$. 
\end{proof}
\begin{proof}[Proof (necessity)]
Suppose, contrary to our claim, that $S^{x}_k$ for some $x \in \rz^N$ or $S^{\infty}_k$ are
finite. Hence there is a bounded sequence $(u_n)_{n \in \nz}$ in $D^{1,2}(\rz^N)$ such that
\begin{align}
\label{eq:38}
\int k |u_n|^q =1, u_n \in D^{1,2}(B_{1/n}(x)) \text{ or } u_n \in D^{1,2}(\rz^N\backslash
B_n(0)).  
\end{align}
Passing to a subsequence we may assume $u_n \schwach u$ and $u_n(x) \to u(x)$ for almost every $x \in
\rz^N$. We conclude from (\ref{eq:38}) that $u \equiv 0$. The compactness of the embedding forces
$\int k(x) |u_n|^q \to 0$ as $n \to \infty$, contrary to (\ref{eq:38}).   
\end{proof}
\begin{remark}
Furthermore (see \cite{SchnDis01}), there are positive constants $c_3(N,q), c_4(N,q)$ such that:
\begin{align*}
&\frac{c_3}{\sqrt{S_k^*}} \le \lim_{\delta \to 0} \supm_{\substack{x \in \rz^N\\0<\rho <\delta}} 
\rho^{(1-\frac{N}{2})} \left(\intg_{B_\rho(x)} k\right)^{\frac{1}{q}} 
\le \frac{c_4}{\sqrt{S_k^*}}  \\
&\frac{c_3}{\sqrt{S^{\infty}_{k}}} \le \lim_{R \to \infty}
\supm_{\substack{x\in \rz^N\\0<\rho}} \rho^{(1-\frac{N}{2})} 
\left(\intg_{B_\rho(x)\backslash B_R(0)} k\right)^{\frac{1}{q}} 
\le \frac{c_4}{\sqrt{S^{\infty}_{k}}} 
\end{align*}
Hence condition (\ref{eq:19}) of Theorem \ref{sec:compactness:t1} prevents point concentration 
whereas condition (\ref{eq:18}) is related to the possible loss of mass at infinity.
\end{remark}
In presence of a group symmetry we denote the length of the orbit
containing $x\in \rz^N$ by
\[|G_x| := \#\{gx\where  g \in G\} \text{ and }|G_\infty|:=1.\]
\begin{corollary}
Suppose $q>2$, $D^{1,2}(\rz^N)$ is continuously embedded in $L^q(\rz^N,k)$ for
some nonnegative $k \in L^1_{loc}$, $G$ is a compact subgroup of $O(N)$
and $k$ is $G$-symmetric. Then $D^{1,2}_G(\rz^N)$ is compactly embedded in $L^{q}(\rz^N,k)$ if 
\[\inf_{x \in \rz^{N}\cup \{\infty\}} \left\{|G_x|\, S^{x}_{k}\right\} = \infty.\]
\end{corollary}
\begin{proof}
Suppose $(u_n)_{n \in \nz}$ is bounded in $D^{1,2}_G(\rz^N)$. We may assume $u_n \schwach 0$ 
in $D^{1,2}_G(\rz^N)$. Because ${\displaystyle D^{1,2}(\rz^N) = D^{1,2}_G(\rz^N) \oplus D_G^{1,2}(\rz^N)^{\bot}}$ 
we have $u_n \schwach 0$ in $D^{1,2}(\rz^N)$. Lemma \ref{sec:conc-comp-lemma-l2}
 yields, that there is an at most countable 
set ${\displaystyle S:= \{x_j\where  j \in J\}}$ of distinct points in ${\displaystyle
\rz^{N}}$ and a family ${\displaystyle \{\nu_j\where  j \in J\}}$ of positive
numbers such that  
\[ k |u_n|^q \schwach \nu = \sum_{j \in J} \nu_j \delta_{x_j} \text{ and } |\nabla u_n|^2 \schwach \mu \ge 
\sum_{j \in J} \mu_j \delta_{x_j}\]
weakly in sense of measures ${\displaystyle S_k^{x_j}\, \nu_j^{2/q} \le \mu_j \text{ and } 
\mu_{\infty} \ge S^{\infty}_{k}\, \nu_\infty^{2/q}}$.\\
Because ${\displaystyle S^{\infty}_{k}= \infty}$ we see $\nu_\infty =0$. Suppose $x_0 \in S \neq \emptyset$.
Then the $G$-symmetry of the involved measures implies ${\displaystyle \{gx_0 \where g \in G\} \subset S}$
and we have
\[ |G_{x_0}| S_k^{x_0} \nu_0^{2/q} \le \sum_{j \in J} \mu_j < \infty.\]
Thus $\nu_0 =0$; a contradiction. Hence Lemma \ref{sec:conc-comp-lemma-l2} $(4)$ leads to 
${\displaystyle \int k |u_n|^q \pfeillang^{n \to \infty} 0}$.
\end{proof}
\section{Palais-Smale condition}
\label{palais}  
In the remainder of this section we always assume (\ref{eq:13})-(\ref{eq:12}).
Because of (\ref{eq:12}) we have ${\displaystyle E_G= D_G^{1,2}(\rz^N)\cap L^{p}(\rz^N,h^{-})}$ 
and we may replace $\|u\|_{E_G}$ 
with the equivalent norm $\|\nabla u\|_2 + \|(h^-)^{1/p} u\|_p$. We still consider the functional
$I: E_G \to \rz$ defined by  
\[I(u) := \frac{1}{2} \|\nabla u\|_2^{2} - \frac{1}{p} \int h(x) |u|^{p}.\]
Clearly ${\displaystyle I \in C^{1}(E_G,\rz)\cap C^{1}(E,\rz)}$ and  
\[I'(u) \phi = \int \nabla u \nabla \phi - \int h(x) |u|^{p-2}u \phi.\]
Critical points of $I$ correspond to weak solutions of
\begin{align}
\label{eq:22}
-\laplace u - h(x) |u|^{p-2} u =0, \quad u \in E_G.  
\end{align}
\begin{lemma}[symmetric criticality]
\label{palais-l1}  
Let $(u_n)_{n \in \nz}$ be a $(PS)_c$ sequence in $E_G$, i.e.  
\[I(u_n) \to c, \qquad I'(u_n) \to 0 \mbox{ in } E_G'.\]  
Then ${\displaystyle I'(u_n) \to 0}$ in $E'$.   
\end{lemma}  
\begin{proof}
The $G$- symmetry of the Laplacian and of $h$ yield the $G$-symmetry of $I$.
Consequently we have for $g \in G$, $u \in E_G$ and $v \in E$ 
\begin{align*}
I'(u) g_* v &= \lim_{t \to 0}\frac{I(u+t g_* v)-I(u)}{t}  = \lim_{t \to 0}\frac{I(g_*^{-1} u+t v)-I(g_*^{-1} u)}{t}\\
&= I'(g_*^{-1}u)v = I'(u)v. 
\end{align*}
Let $\mu$ denote the Haar measure corresponding to the compact group $G$, then we have
\begin{align*}
I'(u)\left(\int_G g_* v \text{d}\mu(g)\right) &= \int_G I'(u) g_* v \text{d}\mu(g) = 
\int_G I'(u) v \text{d}\mu(g)= I'(u) v.
\end{align*}
Hence 
\[ \sup_{\|v\|_E=1} I'(u) v = \sup_{\|v\|_E=1,\, v \in E_G} I'(u)v,\]
and the claim follows.
\end{proof} 

\begin{lemma}[$(PS)_c$ condition] 
\label{palais-l2}  
Every $(PS)_c$ sequence $(u_n)_{n \in \nz}$ in $E_G$ contains a
convergent subsequence if one of the following conditions is satisfied
\begin{align}
\label{eq:31}
\begin{split}
&D^{1,2}_G(\rz^N) \text{ is compactly embedded in }
L^p(\rz^N,h^+);
\end{split}\\[1ex]
\label{eq:32}
\begin{split}
&D^{1,2}(\rz^N) \text{ is continuously embedded in }L^p(\rz^N,h^+)
\text{ and }\\
&c < c_0:=\left(\frac{1}{2}-\frac{1}{p}\right) \inf_{x \in \rz^{N}\cup
\{\infty\}} \left\{|G_x| (S^{x}_{h^{+}})^{\frac{p}{p-2}} \right\}.
\end{split}
\end{align}
\end{lemma}  
\begin{proof}  
Let $(u_n)_{n \in \nz}$ be a sequence in $E_G$ such that ${\displaystyle I(u_n)
\to c}$ and ${\displaystyle I'(u_n) \to 0}$ in ${E_G'}$. Because 
$\|u\|_{E_G}=\|\nabla u_n\|_2+\|(h^{-})^{1/p} u_n\|_p$ we have  
\begin{align}
\label{eq:25}  
c + o(1)\left( \|\nabla u_n\|_2+\|(h^{-})^{1/p} u_n\|_p \right) = 
I(u_n)- (1/p)I'(u_n)u_n= (\frac{1}{2}-\frac{1}{p}) \|\nabla u_n\|_2^{2}.  
\end{align}  
Suppose $\|\nabla u_n\|_2 \to \infty$. Then equation (\ref{eq:25}) implies  
${\|(h^{-})^{1/p} u_n\|_p \ge \|\nabla u_n\|_2^{2}}$  
for large $n$. Hence
\[ c+1 \ge I(u_n) \ge \frac{1}{2} \|\nabla u_n\|_2^{2} + \frac{1}{p} \|\nabla
u_n\|_2^{2p} - C \|\nabla u_n\|_2^{p}\]  
for large $n$, which is impossible. Thus $(u_n)_{n \in \nz}$ is bounded in
$E_G$. Passing to a subsequence we may assume ${\displaystyle u_n \schwach  u}$ in $E_G$.\\
Suppose for a moment there holds 
\begin{align}
\label{eq:33}
u_n \to  u \text{ in } L^{p}(\rz^N,h^+) \text{ as }n \to \infty.  
\end{align}
The fact that $u_n$ converges weakly to $u$ in $E_G$, $D^{1,2}(\rz^N)$ and $L^{p}(\rz^N,h^{-})$ implies
${\displaystyle I'(u)=0}$. Calculating  
\begin{align*}  
0 &= \lim_{n \to \infty} (I'(u_n)-I'(u))(u_n-u)\\  
&= \lim_{n \to \infty} \big(\int |\nabla u_n - \nabla u|^{2} - \int
h^{+}(|u_n|^{p-2}u_n -|u|^{p-2}u)(u_n-u) \\  
& +\int h^{-}\underbrace{(|u_n|^{p-2}u_n -|u|^{p-2}u)(u_n-u)}_{\ge 0}\big)
\\  
&\ge \lim_{n \to \infty} \big(\int |\nabla u_n - \nabla u|^{2} - \int
h^{+}(|u_n|^{p-2}u_n -|u|^{p-2}u)(u_n-u)\big)\\  
&= \lim_{n \to \infty} \int |\nabla u_n - \nabla u|^{2} + o(1)_{n \to
\infty}  
\end{align*}  
we see $u_n \to u$ in $D^{1,2}(\rz^N)$. Finally we have  
\begin{eqnarray*}  
0 &=& \lim_{n \to \infty} I'(u_n)u_n -I'(u)u = \lim_{n \to \infty} \int
h^{-}|u_n|^{p} - \int h^{-}|u|^{p}  
\end{eqnarray*}  
and the uniform convexity of $L^{p}$ implies ${\displaystyle u_n \to u}$ in $E_G$.\\
What is left is to show (\ref{eq:33}), which immediately follows under assumption (\ref{eq:31}).
Thus the proof is completed by showing that (\ref{eq:33}) holds under assumption (\ref{eq:32}).\\
By Lemma \ref{sec:conc-comp-lemma-l2} 
there exist $G$-symmetric measures $\mu$ and $\nu$ satisfying $(1)-(4)$ of
Lemma \ref{sec:conc-comp-lemma-l2} and
\[ |\nabla u_n|^2 \pfeillang^{n \to \infty} \mu, \quad h^+|u_n|^p \pfeillang^{n \to \infty} \nu.\]
Let $x_k$ be an atom of $\nu$. We take
${\displaystyle \phi \in C^{1}(\rz^{N})}$ such that  
\[ \phi \cdot \charak_{B_1(0)} \equiv 1, \quad \phi \cdot (1-\charak_{B_2(0)}) \equiv
0, \quad |\nabla \phi|\le 2\]  
and define ${\displaystyle \phi_\epsilon(x) :=
\phi\left(\frac{x-x_k}{\epsilon}\right)}$. Lemma \ref{palais-l1} implies
${\displaystyle I'(u_n)\phi_\epsilon u_n \to 0}$. Hence
\[\int |\nabla u_n|^2 \phi_\eps + \int \nabla u_n \nabla \phi_\eps u + \int h^- |u_n|^q \phi_\eps
-\int h^+ |u_n|^q \phi_\eps \pfeillang^{n \to \infty} 0.\]
This leads to the following estimate  
\begin{eqnarray*}  
\int \phi_\epsilon d\mu - \int \phi_\epsilon d\nu &\le & \limsup_{n \to
\infty} \int |\nabla u_n| |u_n| |\nabla \phi_\epsilon|\\  
&\le& C \limsup_{n \to \infty} \left(\int |u_n|^{2} |\nabla
\phi_\epsilon|^{2} \right)^{1/2}\\  
&=& C \left(\int |u|^{2} |\nabla \phi_\epsilon|^{2} \right)^{1/2}\\  
&\le& C \|u \cdot \charak_{B_{2\epsilon}(x_k)}\|_{2^{*}}
\underbrace{\left(\int_{B_{2\epsilon}(x_k)} |\nabla \phi_\epsilon|^{N} \right)^{1/N}}_{\le \text{const}}\\  
&\le& o(1)_{\epsilon \to 0}.  
\end{eqnarray*}  
Thus ${\displaystyle \nu_k \ge \mu_k}$ and  Lemma \ref{sec:conc-comp-lemma-l2} (3)
implies:
\begin{align}
\label{eq:26}
\nu_j \ge (S^{x_j}_{h^{+}})^{\frac{p}{p-2}}.  
\end{align}
Take ${\displaystyle \phi_R \in C^{1}(\rz^{N})}$ $G$-symmetric such that  
\[\phi_R(x)=1, \; \forall |x|>R+1 \qquad \phi_R(x)=0, \; \forall |x|<R \qquad
0\le\phi(x) \le 1.\]  
Then  
\begin{eqnarray*}  
0 &=& \lim_{n \to \infty} I'(u_n) \phi_R u_n \\  
&\ge& \limsup_{n \to \infty} \left( \int |\nabla u_n|^{2} \phi_R - \int
|\nabla u_n| |u_n| |\nabla \phi_R| - \int h^{+} |u_n|^{p} \phi_R \right).  
\end{eqnarray*}  
As before we see  
\[\lim_{R \to \infty} \limsup_{n \to \infty} \int |\nabla u_n| |u_n| |\nabla
\phi_R| \le C \lim_{R \to \infty} \|u \cdot \charak_{R<|x|<R+1}\|_{2^{*}} = 0.\]  
Hence ${\displaystyle \nu_\infty \ge \mu_\infty}$ and  
\begin{align}
\label{eq:27}
\nu_\infty \ge (S^{\infty}_{h^-})^{\frac{p}{p-2}}.   
\end{align}
For every $\psi \in C^{\infty}_c(\rz^{N})$ with $0 \le \psi(x) \le 1$ there
holds  
\begin{align}  
c &= \lim_{n \to \infty}\left( I(u_n) - \frac{1}{p} I'(u_n)u_n \right)=
\left(\frac{1}{2}- \frac{1}{p}\right) \lim_{n \to \infty} \int |\nabla
u_n|^{2} \notag \\  
&\ge \left(\frac{1}{2}- \frac{1}{p}\right) \lim_{n \to \infty} \int |\nabla
u_n|^{2} \psi \label{eq:28} 
\end{align}  
Suppose $x_k$ is an atom of $\nu$.  For each $g \in G$, due to the $G$-symmetry, $g x_k$ is 
an atom of $\nu$ with the same mass. The $G$-symmetry of $h$
leads to ${\displaystyle S^{x_k}_{h^{+}}= S^{g x_k}_{h^{+}}}$ for all $g \in
G$. Choose a $\psi$ as above with $\psi(g x_k)=1$ for all $g \in G$. Then (\ref{eq:26}) and (\ref{eq:28})
imply  
\[c \ge \left(\frac{1}{2}- \frac{1}{p}\right) |G_{x_k}| S^{x_k}_{h^{+}}\,
\nu_k^{2/p} \ge \left(\frac{1}{2}- \frac{1}{p}\right) |G_{x_k}|
(S^{x_k}_{h^{+}})^{\frac{p}{p-2}},\]  
which is impossible. Hence ${\displaystyle J = \emptyset}$. If we use
$\phi_R$ in the estimate (\ref{eq:28}), we see with the help of (\ref{eq:27}) 
\[c \ge \left(\frac{1}{2}- \frac{1}{p}\right)
(S^{\infty}_{h^+})^{\frac{p}{p-2}}\]  
and get ${\displaystyle \nu_\infty=0}$. Consequently Lemma \ref{sec:conc-comp-lemma-l2} $(4)$ and
the uniform convexity of $L^{p}$ implies  
\[u_n \to  u \text{ in } L^{p}(\rz^N,h^+) \text{ as }n \to \infty.\]  
\end{proof}
\subsection*{Existence of positive solutions}
We use the mountain pass Lemma of Ambrosetti and Rabinowitz \cite{AmbRab73} and get  
\begin{theorem}
\label{s4t1}
If $D^{1,2}_G(\rz^N) \embed L^p(\rz^N,h^+)$ is compact, then (\ref{eq:10}) has a 
nontrivial, nonnegative weak solution. 
\end{theorem}
\begin{theorem}
\label{s4t2}  
Suppose $D^{1,2}(\rz^N)$ is continuously embedded in $L^p(\rz^N,h^+)$ 
and there is an ${\displaystyle u \in E_G}$ such that  
\[\int h |u|^{p}>0 \mbox{ and } \max_{0 \le t< \infty} I(tu) \le c_0.\]  
Then (\ref{eq:10}) is solvable.
\end{theorem}  
\begin{proof}[Proof of Theorems \ref{s4t1} and \ref{s4t2}]  
$D^{1,2}_G(\rz^N)$ is continuously embedded in $L^p(\rz^N,h^+)$. Consequently
\begin{align}
\label{eq:41}
I(u) = \frac{1}{2} \|\nabla u\|_2^2+\frac{1}{p}\int h^- |u|^p-\frac{1}{p}\int h^+ |u|^p 
\ge \frac{1}{2} \|\nabla u\|_2^2+\frac{1}{p}\int h^- |u|^p - 
C \|\nabla u\|_2^p.
\end{align}
Hence  
\[c := \inf_{\gamma \in \Gamma} \max_{t \in [0,1]} I(\gamma(t)) >0,\]  
where ${\displaystyle \Gamma :=\{\gamma \in C([0,1], E_G)\where   \gamma(0)=0,\,
I(\gamma(1))<0\}}$.\\ 
The mountain pass theorem leads to a $(PS)_c$ sequence.  
If $c<c_0$ or (\ref{eq:31}) holds, then we obtain a critical point $u_0 \in E_G$ of $I$ with the help of
Lemma \ref{palais-l2}.\\  
If $c=c_0$, then the infimum is attained by the path ${\displaystyle
\gamma_0: t \mapsto t t_0 u}$ for a suitable $t_0$. Let $u_0:= \tilde t t_0
u$ with ${\displaystyle I(\gamma_0(\tilde t)) = \max I(\gamma_0(t))}$. Then
${\displaystyle I'(u_0)=0}$, because otherwise $\gamma_0$ can be deformed to
a path $\gamma_1$ with ${\displaystyle \max I(\gamma_1(t))<c}$ contradicting
the definition of $c$.\\  
In both cases we obtain a critical point $u_0$ with $I(u_0)=c>0$. Because all
the involved terms will not change their values, if we replace $u_0$ by
$|u_0|$, we have  
\[c= I(|u_0|) = \max_{t} I(t|u_0|)\]  
and may deduce as above to see, that $|u_0|$ is also a critical point of
$I$.  
\end{proof}  

\begin{corollary} \label{secpaconc1} 
Suppose $D^{1,2}(\rz^N)$ is continuously embedded in $L^p(\rz^N,h^+)$.  
Then (\ref{eq:10}) has a solution if\\ 
${\displaystyle \exists \, u \in E_G:}$ ${\displaystyle \int h
|u|^{p}>0}$, $\|\nabla u\|_2 =1$,  
\begin{align}
\label{eq:7}
\int h |u|^{p} \ge  \sup_{x \in \rz^{N}\cup \{\infty\}} 
\left\{|G_x|^{-\frac{p-2}{2}} (S^{x}_{h^{+}})^{-\frac{p}{2}} \right\}.  
\end{align}
\end{corollary}  

\begin{remark}  
Let $u_0$ be a positive mountain pass solution found in Theorem \ref{s4t1}.
Thanks to the homogeneity of the nonlinear part of $I$ it is easy to see that  
$u_0$ minimizes the functional $I$ on its Nehari manifold, i.e.  
\[I(u_0) = \inf \{ I(u)\where  I'(u)u =0,\; u \neq 0\}.\]  
\end{remark}
 
\subsection*{Existence of multiple solutions}
$I$ is an even functional, i.e. $I(u)=I(-u)$. Consequently classical results leading to multiple
critical points of symmetric functionals apply. We use a version of \cite[Thm. 9.12]{Rab86} given in
\cite{BiaChaSzu95}
\begin{varthm}[Theorem 2 in \cite{BiaChaSzu95}] Let $E$ be an infinite dimensional Banach space and
  $I\in C^1(E,\rz)$ an even functional satisfying $(PS)_c$ condition for each $c$ and
  $I(0)=0$. Furthermore, 
  \begin{itemize}
  \item[(i)] there exists $\alpha,\rho>0$ such that $I(u)\ge \alpha$ for all $\|u\| = \rho$;
  \item[(ii)] there exits an increasing sequence of subspaces $(E_n)_{n \in \nz}$ of $E$, with $\dim
    E_n=n$, such that for every $n$ one can find a constant $R_n>0$ such that $I(u)\le 0$ for all $u
    \in E_n$ with $\|u\| \ge R_n$.
  \end{itemize}
  Then $I$ possesses a sequence of critical values $(c_n)_{n \in \nz}$ tending to $\infty$ as $n \to \infty$. 
\end{varthm}
\begin{theorem}
\label{sec:exist-mult-solut:t1}
Suppose there exits a smooth $G$-symmetric domain $\emptyset \neq \Omega \subset \rz^N$ such that $h(x)>0$ for all $x
\in \Omega$ and $D^{1,2}_G(\rz^N)$ is compactly embedded in $L^p(\rz^N, h^+)$. Then (\ref{eq:22})
has infinitely many $G$-symmetric solutions.
\end{theorem}
\begin{proof}
We apply \cite[Thm. 2]{BiaChaSzu95} with $E=E_G$. From (\ref{eq:41}) we obtain $(i)$. To see that
$(ii)$ holds we choose an increasing sequence  $(E_n)_{n \in \nz}$ of subspaces of
$D^{1,2}_G(\Omega)$ with $\dim E_n = n$. We may assume $D^{1,2}_G(\Omega) \subset E_G$ by extending
the functions by $0$ outside $\Omega$. Since the dimension of each $E_n$ is finite, we conclude
\[\inf_{u \in E_n, \|u\|=1} \, \int h(x) |u|^p =: \epsilon(n)>0,\]
which immediately implies $(ii)$. By Lemma \ref{palais-l2} the $(PS)_c$ condition holds for all $c$,
which completes the proof.    
\end{proof}
\section{Examples}
\label{exa}
It is known (see for instance \cite{ChoChu93,Lie83}) that for $N \ge 3$ and $0\le \delta <2$
\begin{align}
\label{eq:4}
\sup_{u \in D^{1,2}(\rz^N)\backslash \{0\}}
\frac{\left(\int |x|^{-\delta} |u|^p \right)^{\frac{1}{p}}}{\|\nabla u\|_2} = K(N,\delta) < \infty,  
\end{align}
where ${\displaystyle p = p(N,\delta)= \frac{2(N-\delta)}{N-2}}$. $K(N,\delta)$ is attained by the function
\[u_\delta(x) = (1+|x|^{2-\delta})^{-(N-2)/(2-\delta)}.\]
Thanks to the dilatation symmetry of the quotient in (\ref{eq:4}) the functions 
\[u_{\delta, \sigma}(x):= u_\delta(x/\sigma) \text{ and }
v_{\delta,\sigma}:= \frac{u_{\delta,\sigma}}{\|\nabla u_{\delta,\sigma}\|_2}\]
also maximize (\ref{eq:4}) for
all $\sigma>0$.\\
Suppose $k \in C(\rz^N)$ is a continuous nonnegative function. Then we easily get for 
$N \ge 3$, $0\le \delta <2$ and $p=p(N,\delta)$
\begin{align}
\label{eq:5}
S^x_{k(\cdot) |\cdot|^{-\delta}} &= 
\begin{cases}
  \infty& \text{if $\delta >0$ and $x \neq 0$,}\\
  k(x)^{-2/p} K(N,0)^{-2}& \text{if $\delta =0$,}\\
  k(0)^{-2/p} K(N,\delta)^{-2}& \text{if $\delta >0$ and $x=0$.}
\end{cases}
\end{align}
Furthermore, we have 
\begin{align}
\label{eq:6}
\left(\limsup_{|x|\to \infty}k(x)\right)^{-2/p} K(N,\delta)^{-2} \le S^\infty_{k(\cdot) |\cdot|^{-\delta}}    
\le \left(\liminf_{|x|\to \infty}k(x)\right)^{-2/p} K(N,\delta)^{-2}.
\end{align}

\begin{corollary}
\label{sec:examples:c3}
Suppose $h \in L^1_{loc}$, $\int h |u|^p >0$ for some $u \in
C_c^\infty(\rz^N)$.
Then (\ref{eq:10}) has a solution if one of the following conditions holds
\begin{align}
\label{eq:29}
\begin{split}
&2<p<2^*,\, h^+ \le \sum_{i \in \nz} \alpha_i
f(x-p_i)|x-p_i|^{\frac{N-2}{2}p-N} \text{, where } p_i \in \rz^N,\,
\alpha_i \in \rz,\\
&\sum_{i \in \nz}|\alpha_i| <\infty \text{ and } f \in L^\infty(\rz^N)
\text{ such that } f(x) \xpfeil[x \to 0]{|x| \to \infty} 0. 
\end{split}\\
\label{eq:34}
\begin{split}
&p\ge 2^*,\, h \text{ is radially symmetric and satisfies } (\ref{eq:16}).  
\end{split}
\end{align}
\end{corollary}
\begin{proof}
(\ref{eq:34}) implies $D^{1,2}_{O(N)}(\rz^N)$ is compactly embedded in
$L^p(\rz^N,h^+)$ (see for instance \cite[Lem. 6]{Rot90}).\\
Denote by $f_i(x)$ the function $f_i(x):=f(x-p_i)|x-p_i|^{\frac{N-2}{2}p-N}$. 
By Corollary \ref{sec:compactness:co1} the inclusion of $D^{1,2}(\rz^N)$ in $L^p(\rz^N,f_i)$ 
is compact for all $i \in \nz$. Hence we may estimate
\begin{align*}
\supm_{\substack{x \in \rz^N\\0<\rho <\delta}} \rho^{(1-\frac{N}{2})p} 
\int_{B_\rho(x)} \sum_{i \in \nz}\alpha_i f_i
&\le \supm_{\substack{x \in \rz^N\\0<\rho <\delta}} \rho^{(1-\frac{N}{2})p} \int_{B_\rho(x)} f_1 \, 
\left(\sum_{i \in \nz}|\alpha_i|\right) = o(1)_{\delta \to 0},\\
\supm_{\substack{x\in \rz^N\\0<\rho}} \rho^{(1-\frac{N}{2})p} \hspace{-1em}
\intg_{B_\rho(x)\backslash B_R(0)} \sum_{i \in \nz}\alpha_i f_i &\le \supm_{\substack{x\in \rz^N\\0<\rho}} 
\rho^{(1-\frac{N}{2})p} \left(\sum_{i\ge k} |\alpha_i| \int_{B_\rho(x)} f_1  + \hspace{-1em}
\intg_{B_\rho(x)\backslash B_R(0)} \sum_{i<k}\alpha_i f_i\right)\\ 
&\le o(1)_{k \to \infty} + o(1)_{R \to \infty}.
\end{align*}
Thus (\ref{eq:29}) implies $D^{1,2}(\rz^N)$ is compactly embedded in $L^p(\rz^N,h^+)$.\\
In both cases Theorem \ref{s4t1} yields the existence of a solution.
\end{proof}
Pohozaev's identity adapted to (\ref{eq:22}) leads to
\begin{lemma}
\label{sec:examples:l1}
Suppose $2<p$ and $h(x)= k(x)|x|^{-\delta}$ for some $k \in
C^1(\rz^N)$, where $\delta = N -\frac{p}{2}(N-2)$. Then every solution
$u \in C^2(\rz^N)$ of (\ref{eq:22}), such that the function $<\nabla k(x),x>\,|x|^{-\delta} |u(x)|^p
\in L^1(\rz^N)$, satisfies 
\begin{align}
\label{eq:8}
\int <\nabla k(x),x>|x|^{-\delta} |u|^p =0.   
\end{align}
\end{lemma}
\begin{proof}
Suppose $u \in D^{1,2}(\rz^N) \cap L^p(|h|) \cap C^2(\rz^N)$ solves (\ref{eq:22}). Because 
\begin{align}
\label{eq:35}
<\nabla h(x),x> = <\nabla k(x),x>|x|^{-\delta}- \delta h(x)  
\end{align}
we have $\left|<\nabla h(x),x>\right| \,|u|^p \in L^1(\rz^N)$. We use a version of Pohozaev's
identity \cite{Poh70} given in \cite[Thm. 29.4]{KuzPoh97} and (\ref{eq:35}) to derive 
\begin{align}
\label{eq:36}
\begin{split}
\frac{N-2}{2} \int |\nabla u|^2 &= \frac{N}{p} \int h(x) |u|^p + \frac{1}{p} \int <\nabla h(x),x>
|u|^p\\
 &= \frac{N-2}{2} \int h(x) |u|^p + \frac{1}{p} \int <\nabla k(x),x> |u|^p. 
\end{split}
\end{align}
The fact that $I'(u)u =0$ and (\ref{eq:36}) give the identity (\ref{eq:8}).
\end{proof}

\begin{corollary}
\label{sec:examples:c4}
Consider the equation
\begin{align}
\label{eq:37}
-\laplace u = (1+|x|)^{-\delta} |u|^{p-2}u, \quad 0\neq u \in D^{1,2}(\rz^N) \cap
 L^p(\rz^N,(1+|x|)^{-\delta}).
\end{align}
\begin{itemize}
\item[(i)] (\ref{eq:37}) has no solution $u \in C^2(\rz^N)$ if $2<p<2^*$ and $\delta \le N
  -\frac{p}{2}(N-2)$.
\item[(ii)] (\ref{eq:37}) has infinitely many $C^2(\rz^N)$-solutions if $2<p<2^*$ and $\delta > N
  -\frac{p}{2}(N-2)$. At least one solution is strictly positive in $\rz^N$.  
\end{itemize}
\end{corollary}
\begin{proof}
From $2<p<2^*$, $\delta > \delta_0:= N  -\frac{p}{2}(N-2)$ and Corollary \ref{sec:examples:c3} we conclude that
(\ref{eq:37}) possesses a nontrivial weak solution $u$. Due to Harnack's inequality and standard
regularity results (see \cite[C.3]{Sim82}) $u$ is a strictly positive
$C^2(\rz^N)$-function. In addition by Theorem \ref{sec:exist-mult-solut:t1} there are infinitely
many solutions of (\ref{eq:37}).\\ 
Considering $k(x):= (1+|x|)^{-\delta} |x|^{\delta_0}$ it is easy to check that Lemma
\ref{sec:examples:l1} yields the desired nonexistence result.   
\end{proof}
Solutions of (\ref{eq:10}) may be obtained in a non-compact setting if it is possible to find appropriate
test functions to ensure that the mountain pass value $c$ is below the compactness threshold $c_0$, i.e. we 
have to find a function $u \in E_G$ that satisfies (\ref{eq:7}). 
\begin{remark}
\label{sec:examples:r1}
Suppose $2<p\le 2^*$ and $h \in L^1_{loc}$ is $G$-symmetric such that 
\[h^+(x) = k(x) |x|^{-\delta},\]
where $k \in C(\rz^N)$ and ${\displaystyle \delta = N -\frac{p}{2}(N-2)}$. With the notation
\[k_c := \max(\sup_{x \in \rz^N}\{k(x)|G_x|^{-\frac{p-2}{2}}\}, \limsup_{|x|\to \infty}k(x))\]
we have
\begin{align*}
\sup_{x \in \rz^{N}\cup \{\infty\}} 
\left\{|G_x|^{-\frac{p-2}{2}} (S^{x}_{h^{+}})^{-\frac{p}{2}} \right\} \le k_c K(N,\delta)^p= 
k_c \int |x|^{-\delta} v_{\delta,\sigma}^p.  
\end{align*}
Thus (\ref{eq:7}) holds if we have for some $\sigma >0$
\begin{align}
\label{eq:9}
\int h(x)v_{\delta,\sigma}^p >0,\; \int (h(x)-k_c|x|^{-\delta})v_{\delta,\sigma}^p \ge 0.  
\end{align}
\end{remark}
Under the assumptions of Remark \ref{sec:examples:r1} the space $D^{1,2}_G(\rz^N)$ is not compactly embedded in 
$L^p(\rz^N,h^+)$ if $k(0)>0$ or $\liminf_{|x| \to \infty} k(x)>0$. Corollaries \ref{sec:examples:c1} and 
\ref{sec:examples:c2} below yield some sufficient conditions for the existence of solutions to (\ref{eq:10})
in the non-compact case. We leave it to the reader to verify with the help of (\ref{eq:9}) that the
proofs given in \cite[Cor. 1,2]{BiaChaSzu95} carry over to our situation.  
\begin{corollary}
\label{sec:examples:c1}
Assuming the hypotheses of Remark \ref{sec:examples:r1} and
\begin{align*}
k(0) \ge 
\begin{cases}
\limsup_{|x|\to \infty}k(x)& \text{if } 2<p<2^*,\\
\max\{\limsup_{|x|\to \infty}k(x),\sup_{x \in \rz^N}|G_x|^{-\frac{p-2}{2}} k(x)\}& \text{if } p=2^*  
\end{cases}
\end{align*}
(\ref{eq:10}) has a solution if one of the following conditions is satisfied
\begin{align*}
&h\neq 0, \,h(x) \ge k(0) |x|^{-\delta} \text{ for all }x \in \rz^N\backslash \{0\};\\
&\exists \epsilon,r >0:\: \int_{|x|\ge r} h^- |x|^{2(\delta -N)} < \infty,\,k(x)\ge k(0) + \eps |x|^{N-\delta}\; 
\forall x \in B_r(0);\\
&\int \left|h(x)-k(0)|x|^{-\delta}\right| \,|x|^{2(\delta-N)} < \infty \text{ and }  
\int (h(x)-k(0)|x|^{-\delta}) |x|^{2(\delta-N)}>0. 
\end{align*}
\end{corollary}
\begin{corollary}
\label{sec:examples:c2} 
Assuming the hypotheses of Remark \ref{sec:examples:r1} and 
\begin{align*}
\lim_{|x|\to \infty}k(x)=: k(\infty) \ge 
\begin{cases}
k(0)& \text{if } 2<p<2^*,\\
\sup_{x \in \rz^N}\{|G_x|^{-\frac{p-2}{2}} k(x)\}& \text{if } p=2^*  
\end{cases}
\end{align*}
(\ref{eq:10}) has a solution if one of the following conditions is satisfied
\begin{align*}
&\exists \epsilon,R >0:\: \int_{|x|\le R} h^- < \infty,\,k(x)\ge k(\infty) + \eps |x|^{-N+\delta}\; 
\forall x \in \rz^N\backslash B_R(0);\\
&\int \left|h(x)-k(\infty)|x|^{-\delta}\right| < \infty \text{ and }  
\int (h(x)-k(\infty)|x|^{-\delta}) >0. 
\end{align*} 
\end{corollary}

\section*{Acknowledgments}
I would like to thank my thesis advisor, Hans-Peter Heinz, for many helpful suggestions during the
preparation of the paper.

\bibliographystyle{amsplain}
\bibliography{comp_indef}

\end{document}